\documentclass{amsproc}
\usepackage{graphicx}
\usepackage{amssymb}
\usepackage{epstopdf}
\usepackage{amsmath,amsthm,amscd,amssymb}
\usepackage{latexsym}
\usepackage[colorlinks,citecolor=red,pagebackref,hypertexnames=false]{hyperref}
\usepackage{geometry}                
\numberwithin{equation}{section}

\theoremstyle{plain}
\newtheorem{theorem}{Theorem}[section]
\newtheorem{lemma}[theorem]{Lemma}
\newtheorem{corollary}[theorem]{Corollary}
\newtheorem{proposition}[theorem]{Proposition}

\theoremstyle{definition}
\newtheorem{definition}[theorem]{Definition}

\newtheorem{case[theorem]}{Case}

\theoremstyle{remark}
\newtheorem{remark}[theorem]{Remark}

\numberwithin{equation}{section}

\begin{document}

\title{\parbox{14cm}{\centering{Fourier integral operators, fractal sets, and the regular value theorem}}}


\author{Suresh Eswarathasan, Alex Iosevich and Krystal Taylor}

\date{today}

\email{suresh@math.rochester.edu} 
\email{iosevich@math.rochester.edu}
\email{taylor@math.rochester.edu}
\address{Department of Mathematics, University of Rochester, Rochester, NY}

\thanks{The work of the second listed author was partially supported by the NSF Grant DMS10-45404}

\begin{abstract} We prove that if ${\mathcal E} \subset {\Bbb R}^{2d}$, for $d \ge 2$, is an Ahlfors-David regular product set of sufficiently large Hausdorff dimension, denoted by $dim_{{\mathcal H}}({\mathcal E})$, and $\phi$ is a sufficiently regular function, then the upper Minkowski dimension of the set 
$$ \{w \in {\mathcal E}: \phi_l(w)=t_l; 1 \leq l \leq m \}$$ does not exceed $dim_{{\mathcal H}}({\mathcal E})-m$, in line with the regular value theorem from the elementary differential geometry. Our arguments are based on the mapping properties of the underlying Fourier integral operators and are intimately connected with the Falconer distance conjecture in geometric measure theory. We shall see that our results are, in general, sharp in the sense that if the Hausdorff dimension is smaller than a certain threshold, then the dimensional inequality fails in a quantifiable way. The constructions used to demonstrate this are based on the distribution of lattice points on convex surfaces and have connections with combinatorial geometry. 
\end{abstract} 

\maketitle


\section{Introduction}

\vskip.125in 

The regular value theorem in elementary differential geometry says that if $\phi: X \to Y$, where $X$ is a smooth manifold of dimension $n$ and $Y$ is a smooth manifold of dimension $m<n$ with $\phi$ a submersion on the set 
\begin{equation}
\{x \in X: \vec{\phi}(x)=y \}, 
\end{equation} for $y$ is a fixed element of $Y$, then the set 
\begin{equation}
{\vec{\phi}}^{-1}(y)=\{x \in X: \vec{\phi}(x)=y \}
\end{equation} is either empty or is a $n-m$ dimensional submanifold of $X$. 

In this paper we consider the situation where $Y={\Bbb R}^m$ and $X$ is replaced by $E \times E$, where $E \subset {\Bbb R}^d$ is a set of a given Hausdorff dimension, which, in general, is far from being a smooth manifold. A direct analog of the regular value theorem would be a statement that the set 
$$\{(x,y) \in E \times E: \phi_l(x,y)=t_l; 1 \leq l \leq m \}$$ is either empty or has fractal dimension exactly $2s-m$, where $s$ is the Hausdorff dimension of $E$. We are able to show, under some reasonable hypotheses on $\phi$, that the upper Minkowski dimension of 
$$\{(x,y) \in E \times E: \phi_l(x,y)=t_l, 1 \leq l \leq m \}$$ does not exceed $2s-m$. 

Generalizations of the regular value theorem have been recently considered in the context of Banach spaces. See \cite{AN09} and the references contained therein. 

To put these ideas into context, we recall that Falconer \cite{Fal86} formulated the now celebrated Falconer distance conjecture, which says that if the Hausdorff dimension of a compact set $E$ in ${\Bbb R}^d$, for $d \ge 2$, is greater than $\frac{d}{2}$, then the Lebesgue measure of the distance set 
$$\Delta(E)=\{|x-y|: x,y \in E\}$$ is positive. See also \cite{Mat85}, \cite{Mat87} and \cite{M95} for related results.  In \cite{Fal86}, Falconer proved that the conclusion holds if the Hausdorff dimension of $E$ is greater than $\frac{d+1}{2}$ by showing that for $\mu$, a probability measure on $E$, 
\begin{equation} \label{falconer} 
\mu \times \mu \{(x,y): t \leq |x-y| \leq t+\epsilon \} \lesssim \epsilon.
\end{equation} 

The key to (\ref{falconer}), though Falconer did not express himself using this language, is the $L^2({\Bbb R}^d) \to L^2_{\frac{d-1}{2}}({\Bbb R}^d)$ bound for the averaging operator 
$$ Tf(x)=\int f(x-y) d\sigma_t(y),$$ where $d\sigma_t$ is the Lebesgue measure on the sphere of radius $t$ and $L^2_s({\Bbb R}^d)$ denotes the usual $L^2$-Sobolev space of $L^2$ functions with $s$ generalized derivatives in $L^2({\Bbb R}^d)$. See \cite{So93} for Sobolev estimates for geometric averaging operators. 


In this paper we shall see that under some reasonable assumptions on smooth functions $\phi_l: {\Bbb R}^d \times {\Bbb R}^d \to {\Bbb R}$, a suitable analog of (\ref{falconer}) can be used to prove the geometric inequality 
\begin{equation} \label{falconerpairedgeneral} \overline{dim}_{{\mathcal M}} \{(x,y) \in E \times E: \phi_l(x,y)=t_l; 1 \leq l \leq m \} \leq 2 \cdot dim_{{\mathcal H}}(E)-m, \end{equation} provided that the Hausdorff dimension of $E$ is sufficiently large in the sense to be quantified below. We note that in the model case when $m=1$ and $\phi(x,y)=|x-y|$, this connection is explored in \cite{CEHIT10}. Before formulating the results, we introduce the main analytic tool used in this paper, the generalized Radon transform. 

\subsection{Generalized Radon transforms} 

Given $f: {\Bbb R}^d \to {\Bbb R}$, define 
\begin{equation} \label{keyoperator} T_{\vec{\phi}_t}f(x):=\int_{\{\phi_l(x,y)=t_l; 1 \leq l \leq m\}} f(y)\psi(x,y) d\sigma_{x,t}(y), \end{equation} where $d\sigma_{x,t}$ is the Lebesgue measure on the set $\{y: \phi_l(x,y)=t_l; 1 \leq l \leq m\}$ and $\psi$ is a smooth cut-off function.  Here $\vec{\phi}=(\phi_1, \dots, \phi_m)$ and $t=(t_1, \dots, t_m)$. We shall assume throughout the rest of this paper
\begin{equation} \label{fibration} \{(\nabla_x \phi_l(x,y))\}_{l=1}^m \text{ and } \{\nabla_y \phi_l(x,y))\}_{l=1}^m \end{equation} form two linearly independent sets of vectors in $\mathbb{R}^{d}$ in a neighborhood of the sets 
\begin{equation}
\{x: \phi_l(x,y)=t_l; 1 \leq l \leq m\} \text{ and } \{y: \phi_l(x,y)=t_l; 1 \leq l \leq m\},
\end{equation} respectively. This can be justified by details in the note of Phong and Stein \cite{PhSt86} and is meant to provided an underlying smooth structure.  We call $T_{\vec{\phi}_t}$ the Radon transform associated to $\vec{\phi}$.  More precisely, 
\begin{equation}
T_{\vec{\phi}_t}f(x): C^{\infty}(\mathbb{R}_y^n) \rightarrow C^{\infty}(\mathbb{R}_x^n \times \mathbb{R}^m_t).\end{equation}  For the purposes of this paper, we treat $t$ as fixed a parameter.  The article \cite{PhSt86} treats these operators in more generality and provides their basic theory.  

\subsection{Main results} 

Given $E$ compact in ${\Bbb R}^d$, for $d \ge 2$, define
\begin{equation} \label{paired} S_t^{\vec{\phi}}(E)=\{(x,y) \in E \times E: \phi_l(x,y)=t_l; 1 \leq l \leq m \}.\end{equation}

Recall that $E \subset {\Bbb R}^d$ is said to be Ahlfors-David regular if there exists a Borel measure $\mu$, supported on $E$ and $C>0$, such that for all $x \in E$,
\begin{equation} \label{adregular} C^{-1} \delta^s \leq \mu(B_{\delta}(x)) \leq C \delta^s \end{equation} for every $\delta>0$, where $s$ is the Hausdorff dimension of $E$ and $B_{\delta}(x)$ is the ball of radius $\delta$ centered at $x$. 

\begin{theorem} \label{main} 
Let $E \subset {\Bbb R}^d$, for $d \ge 2$, be compact and Ahlfors-David regular. Choose a smooth $\vec{\phi}$ such that $T_{\vec{\phi}_t}$ and $S_t^{\vec{\phi}}(E)$ are as in (\ref{keyoperator}) and (\ref{paired}), respectively. Suppose 
\begin{equation} \label{sobolev} T_{\vec{\phi}_t}: L^2({\Bbb R}^d) \to L^2_s({\Bbb R}^d) \end{equation} with constants uniform in $t \in T=T_1 \times T_2 \times \dots \times T_m$, $T_j$ an interval in ${\Bbb R}$, for some $s>0$ and assume 
\begin{equation} \label{dimensionrestriction} 
dim_{{\mathcal H}}(E)>d-s. 
\end{equation}  Then  for $t \in T$, 
\begin{equation} \label{dimensioninequality} \overline{dim}_{{\mathcal M}}(S_t^{\vec{\phi}}(E)) \leq 2 \cdot \dim_{{\mathcal H}}(E)-m. \end{equation} 
\end{theorem} 

\begin{remark} Our method easily extends to the situation where $E \times E$ is replaced by $E \times F$, with the right hand side of (\ref{dimensioninequality}) replaced by $dim_{{\mathcal H}}(E)+dim_{{\mathcal H}}(F)-m$. It is also not particularly essential for our method that both $E$ and $F$ are subsets of the same Euclidean space ${\Bbb R}^d$. We can take $E \subset {\Bbb R}^{d_1}$ and $F \subset {\Bbb R}^{d_2}$. However, due to the current state of knowledge of Sobolev bounds for generalized Radon transforms, our best results are in the case when $d_1=d_2$, making our hypotheses reasonable. \end{remark} 

\begin{remark} It would be very interesting to extend our result to sets of the form   
$$ \{w \in {\mathcal E}: \vec{\phi}(w)=\vec{t} \},$$ where ${\mathcal E} \subset {\Bbb R}^n$ is of a sufficiently large Hausdorff dimension and $\vec{\phi}$ is sufficiently regular. This requires a rather intricate analysis of the Hausdorff dimension of projections of ${\mathcal E}$ and tensor product properties of the resulting measures. We hope to address this issue in a sequel.  \end{remark} 

The following definition is stated in \cite{PhSt86}.
\begin{definition} We say that $\phi: {\Bbb R}^d \times {\Bbb R}^d \to {\Bbb R}$ satisfies the Phong-Stein rotational curvature condition at $t$ if 
\begin{equation} \label{curvecond}
det \begin{pmatrix} 
 0 & \nabla_{x}\phi \\
 -{(\nabla_{y}\phi)}^{T} & \frac{\partial^2 \phi}{dx_i dy_j}
\end{pmatrix}
\neq 0
\end{equation} on the set $\{(x,y): \phi(x,y)=t \}$.
\end{definition} 

We now list some corollaries of Theorem \ref{main} designed to illustrate concrete situations where the degree of smoothing of the operator $T_{\vec{\phi}_t}$ can be explicitly calculated. This list is not meant to be exhaustive but to simply illustrate the range of applicability of our methods. 

\begin{corollary} \label{phongstein} Suppose that $m=1$ and $\vec{\phi} = \phi: \mathbb{R}^d \rightarrow \mathbb{R}$ satisfies the Phong-Stein rotational curvature condition. Then the conclusion of Theorem \ref{main} holds under the assumption that the Hausdorff dimension of $E$ is greater than $\frac{d+1}{2}$. \end{corollary} 

\begin{corollary} \label{curves} Suppose that $m=d-1$ and 
\begin{equation}
\phi_l(x,y)=(x_{l+1}-y_{l+1})-\gamma_l(x_1-y_1),
\end{equation} where the curve 
\begin{equation}
\Gamma=\{(s, \gamma_1(s), \dots, \gamma_{d-1}(s)): s \in [0,1] \}
\end{equation} has non-vanishing curvature and torsion.  Then (\ref{dimensioninequality}) holds if the Hausdorff dimension of $E$ is greater than $d-\frac{1}{d}$.
\end{corollary} 

We shall prove Corollary \ref{phongstein} below. To prove Corollary \ref{curves}, observe that by the van der Corput Lemma (\cite{St93, So93}) if $\sigma_{\Gamma}$ denotes the Lebesgue measure on $\Gamma$, then 
\begin{equation}
|\widehat{\sigma}_{\Gamma}(\xi)| \lesssim {|\xi|}^{-\frac{1}{d}}.
\end{equation}  It follows that (\ref{sobolev}) holds with $s=\frac{1}{d}$ and thus Corollary \ref{curves} follows from Therem \ref{main}.  

We are able to consider more general families of curves under a variety of geometric assumptions. The reference \cite{GSe04} and those contained therein give a thorough description of such estimates. 

\begin{remark} It would be very interesting to consider the set 
$$ \{(x^1, \dots, x^k) \in E_1 \times E_2 \times \dots \times E_k: \phi_l(x^1, \dots, x^k)=t_l, l=1, \dots, m \}$$ and prove that the upper Minkowski dimension of this does not exceed 
$$ dim_{{\mathcal H}}(E_1)+\dots+dim_{{\mathcal H}}(E_k)-m.$$ 

An natural approach to this question, in view of this paper, is via regularity properties of multi-linear variants of generalized Radon transforms. These are operators of the form 
$$ M(f_1, \dots, f_{k-1})(x)=\int \dots \int f_1(y^1) \dots f_{k-1}(y^{k-1}) dK(x,y^1, \dots, y^{k-1}),$$ where $dK$ is a smooth cut-off function times the Lebesgue measure on the set 
$$ \{(x, y^1, \dots, y^{k-1}): \phi_l(x, y^1, \dots y^{k-1})=t_l; 1 \leq l \leq m\}.$$ 

Some special cases of these operators have been studied in conjunction with the study of finite point configuration and Falconer type problems; see \cite{GI10, CEHIT10}. However, nothing resembling a general theory of such operators is currently available. We hope to address this issue in a subsequent paper. \end{remark} 

\subsection{Applications to the variable coefficient Falconer distance problem} Again, the Falconer distance conjecture states that if the Hausdorff dimension of $E$ is greater than $\frac{d}{2}$, then the Lebesgue measure of the set of distances, $\Delta(E)=\{|x-y|: x,	y \in E\}$ is positive. The best known results, due to Wolff \cite{W99} for $d=2$ and Erdogan \cite{Erd05} for $d > 2$, say that if the Hausdorff dimension of $E$ is greater than $\frac{d}{2}+\frac{1}{3}$, then the Lebesgue measure of $\Delta(E)$ is positive. Techniques of this paper allow us to extend Falconer's result to a variable coefficient setting. More precisely, let $\phi: {\Bbb R}^d \times {\Bbb R}^d \to {\Bbb R}$ be a metric on ${\Bbb R}^d$ and define 
$$ \Delta_{\phi}(E)=\{\phi(x,y): x,y \in E \}.$$ 

The main result of this subsection is the following: 
\begin{theorem} \label{variablefalconer} Suppose that $\phi: {\Bbb R}^d \times {\Bbb R}^d \to {\Bbb R}$ is a metric on ${\Bbb R}^d$ satisfying the rotational curvature condition of Phong and Stein described above. Let $E$ be a compact subset of ${\Bbb R}^d$, $d \ge 2$, of Hausdorff dimension greater than $\frac{d+1}{2}$. Then the Lebesgue measure of $\Delta_{\phi}(E)$ is positive. 
\end{theorem} 

\begin{remark} Theorem \ref{variablefalconer} opens the door to a systematic study of the Falconer distance problem on Riemannian manifolds. This has already led us to some interesting connection with the sharp Weyl formula (see e.g. \cite{So93}). We shall address this issue in a subsequent paper (\cite{EIT10.II}). \end{remark} 

To prove Theorem \ref{variablefalconer}, observe that the proofs of Theorem \ref{main} and Corollary \ref{phongstein} above imply that if the Hausdorff dimension of $E$ is greater than $\frac{d+1}{2}$ and $\mu$ is a Frostman measure on $E$, then 
$$ \mu \times \mu \{(x,y) \in E \times E: t-\epsilon \leq \phi(x,y) \leq t+\epsilon \} \lesssim \epsilon.$$ 

Now, for any cover of $\Delta_{\phi}(E)$ by intervals $(t_j, t_j+\epsilon_j)$, 
\begin{equation}
 1=\mu \times \mu(E \times E) \leq \sum_j \mu \times \mu \{(x,y): t_j \leq \phi(x,y) \leq t_j+\epsilon_j \} \lesssim \sum_j \epsilon_j
 \end{equation} and it follows that there exists a uniform constant $c$ such that 
\begin{equation}
 \sum_j \epsilon_j \gtrsim c>0
 \end{equation} for any covering of $\Delta_{\phi}(E)$. Thus the Lebesgue measure of $\Delta_{\phi}(E)$ is positive. 

\subsection{Sharpness of results} There are at least two notions of sharpness that could be discussed in this context. The first and the most important question is to find a threshold $\alpha_0$ such that if the Hausdorff dimension of $E$ is smaller than $\alpha_0$, then the conclusion of Theorem \ref{main} does not in general hold.  A result of this type is proved in Section \ref{finalsection} but is stated below. In the same section we shall discuss the extent to which our main technical estimate (\ref{mama}) is best possible. 

\begin{theorem} \label{sharp} Let $m=1$. There exists $\phi: {\Bbb R}^d \times {\Bbb R}^d \to {\Bbb R}$ satisfying (\ref{phongstein}) with the following property: for every $s<\frac{d+1}{2}$ there exists a set $E$ of Hausdorff dimension $s$ such that 
\begin{equation} \label{bad} \overline{dim}_{{\mathcal M}}(S_t^{\phi}(E))>2 \cdot dim_{{\mathcal H}}(E)-1. \end{equation} 
\end{theorem} 

As the reader shall see below, we set $\phi(x,y)={||x-y||}_B$, where $B$ is a paraboloid and ${|| \cdot ||}_B$ denotes the norm induced by $B$. It is important to note that the best sharpness example we are able to construct for the function $\phi(x,y)=|x-y|$ only show that $(\ref{bad})$ holds if $s<\frac{d}{2}$, instead of $s<\frac{d+1}{2}$ in the case of the paraboloid induced metric. We do not know whether this is merely an artifact of our method, or whether there is indeed a distinction between the Euclidean metric and the metric induced by the paraboloid $B$. A related construction can be found in \cite{BBCRV07} in the context of Fourier averages. 

This shows that Corollary \ref{phongstein} cannot, in general, be improved. The construction used to obtain Theorem \ref{sharp} can be extended to treat the case of $m>1$. 

\subsection{Organization} This paper is organized as follows. In Section \ref{mainsection} below, we prove Theorem \ref{main}. In Section \ref{corollarysection}, we establish Corollary \ref{phongstein}. In the final part of the paper, Section \ref{finalsection}, we discuss the extent to which our results are optimal. 

\vskip.25in 

\section{Proof of the main result} 
\label{mainsection} 

\begin{lemma} \label{bigreduction} 
Let $E \subset \mathbb{R}^d$ be a compact Ahlfors-David regular set of Hausdorff dimension $\alpha > 0$ and $\mu$ be a Frostman measure on $E$.  If 
\begin{equation} \label{mama} \mu \times \mu \{(x,y) \in E \times E: t_l \leq \phi_l(x,y) \leq t_l + \varepsilon, \text{ } 1 \leq l \leq m \} \lesssim \varepsilon^m,\end{equation} then the conclusion of Theorem \ref{main} holds. \end{lemma}

To prove the lemma note that since each $\phi_l$ is Lipschitz, the $\epsilon$ neighborhood of $S_{\vec{t}}^{\vec{\phi}}(E)$, denoted by ${\left(S_{\vec{t}}^{\vec{\phi}}(E)\right)}^{\epsilon}$, is contained in the set 
\begin{equation}
 \{(x,y) \in E^{\epsilon} \times E^{\epsilon}: t_l \leq \phi_l(x,y) \leq t_l + \varepsilon, 1 \leq l \leq m \},
 \end{equation}
 where $E^{\epsilon}$ denotes the $\epsilon$-neighborhood of $E$, and thus 
\begin{equation}
 \mu \times \mu \left\{ {\left(S_{\vec{t}}^{\vec{\phi}}(E)\right)}^{\epsilon} \right\} \leq \mu \times \mu \{(x,y) \in E \times E: t_l \leq \phi_l(x,y) \leq t_l + \varepsilon, 1 \leq l \leq m \} \lesssim \epsilon^m, 
 \end{equation}
 where the last inequality follows from (\ref{mama}) and $\mu$ being supported only on the set $E$. 

On the other hand, for $\epsilon$ sufficiently small, (\ref{adregular}) implies that 
$$ \mu \times \mu \left\{ {\left(S_{\vec{t}}^{\vec{\phi}}(E)\right)}^{\epsilon} \right\} \gtrapprox \epsilon^{2\alpha} \epsilon^{-\gamma},$$ where $\gamma$ is the upper Minkowski dimension of 
$S_{\vec{t}}^{\vec{\phi}}(E)$. We conclude that 
$$ \overline{dim}_{{\mathcal M}}(S_{\vec{t}}^{\vec{\phi}}(E)) \leq 2\alpha-m,$$ as desired. 

Hence, the proof of Theorem \ref{main} reduces to the following claim. 

\begin{proposition}\label{mamaestprop}
Let $\mu$ be a probability measure on $E$ and $T_{\vec{\phi}_t}: L^2({\Bbb R}^d) \to L^2_s({\Bbb R}^d)$ with $d-s < \alpha < d$, for $\alpha=dim_{{\mathcal H}}(E)$, the Hausdorff dimension of $E$. Then 
\begin{equation} \label{2.2}
\mu \times \mu \{(x,y) \in E \times E: t_l \leq \phi_l(x,y) \leq t_l + \varepsilon, 1 \leq l \leq m \} \lesssim \varepsilon^m.
\end{equation}
\end{proposition}

To prove this, take Schwartz the class functions $\eta_0(\xi)$ supported in the ball $\{ |\xi| \leq 4 \}$ and $\eta(\xi)$ supported in the annulus 
\begin{equation}
\{ 1 < |\xi| < 4\} \ \text{with} \ \eta_j(\xi)=\eta(2^{-j}\xi) \text{ for } j \geq 1
\end{equation}
with
\begin{equation}
\eta_0(\xi) + \sum_{j=1}^{\infty} \eta_j(\xi) =1.
\end{equation}  Define the Littlewood-Paley piece of $\mu_j$ by the relation 
\begin{equation}
\widehat{\mu}_j(\xi)=\widehat{\mu}(\xi) \eta_j(\xi).
\end{equation}

Consider the left hand side of (\ref{2.2}).  This can be rewritten as 

\begin{equation} \label{2.3}
\sum_{j,k} \int \int_{\{t_l \leq \phi_l(x,y) \leq t_l + \varepsilon: 1 \leq l \leq m \}} \psi(x,y) d\mu_j(x) \hskip.1cm d\mu_k(y)= \sum_{j,k} \langle \mu_j,T^{\varepsilon} \mu_k \rangle
\end{equation}

\noindent where $<\cdot,\cdot>$ denotes the $L^2(\mathbb{R}^d)$ inner product and 

\begin{eqnarray} \label{expanded}
\nonumber T^{\varepsilon}\mu_k(x)&=&\int_{\{t_l \leq \phi_l(x,y) \leq t_l + \varepsilon: 1 \leq l \leq m \}}\psi(x,y) d\mu_k(y) \\
&=& \int_{t_1}^{t_1 + \varepsilon} ... \int_{t_m}^{t_m + \varepsilon}  \int_{\vec{\phi}(x,y)=r}\psi(x,y) \mu_k(y) d\sigma_{x,r}(y) dr_1 ... dr_m, 
\end{eqnarray} where $d\sigma_{x,r}$ is the Lebesgue measure on the set $\{y: 
\vec{\phi}(x,y)=r\}$ and $r = (r_1,...,r_m)$. It should be noted that the innermost integral on the right side of (\ref{expanded}) is just $T_{\vec{\phi}_r}$ applied to $\mu_k$. It follows that the right hand side of (\ref{2.3}) becomes
$$\sum_{j,k} \int_{t_1}^{t_1+\epsilon} \cdots \int_{t_m}^{t_m+\epsilon}\langle \mu_j, T_{\vec{\phi_r}}(\mu_k)\rangle dr_1 \cdots dr_m.$$  

We will now use the mapping properties of $T_{\vec{\phi_r}}$ to prove $<\mu,T_{\vec{\phi_r}} \mu>$ is uniformly bounded in $r$ over the the domain of integration.  This, in turn, will prove our desired theorem.

We have 

\begin{equation} \label{2.5}
\hskip.1cm \langle \mu, T_{\vec{\phi_r}}(\mu)\rangle \hskip.1cm = \sum_{j,k} \langle \mu_j, T_{\vec{\phi_r}}(\mu_k)\rangle 
\end{equation}

\begin{equation} \label{2.6}
= \hskip.1cm \sum_{|j-k| \leq K} \langle \mu_j, T_{\vec{\phi_r}}(\mu_k)\rangle + \sum_{|j-k| > K} \langle \mu_j, T_{\vec{\phi_r}}(\mu_k)\rangle 
\end{equation}

\noindent for $K$ large enough; the choice of $K$ will be justified later.  We will estimate each of the above sums separately.  For the first sum, 

\begin{equation} \label{mickey} \sum_{|j-k| \leq K}\langle \mu_j, T_{\vec{\phi_r}}(\mu_k)\rangle  \hskip.1cm \lesssim \sum_{|j-k| \leq K}2^{j\frac{d-\alpha}{2}}2^{k\frac{d-\alpha}{2}}2^{-ks}\lesssim 1  \end{equation}
provided that $d-s < \alpha < d$.  Indeed, as $\eta_j \sim \eta_j^2$,

\begin{eqnarray}
\nonumber \sum_{|j-k| \leq K}\langle \mu_j, T_{\vec{\phi_r}}(\mu_k)\rangle &=& \sum_{|j-k| \leq K}\langle \widehat{\mu_j}, \widehat{T_{\vec{\phi_r}}(\mu_k)}\rangle\\
\nonumber &\sim&  \sum_{|j-k| \leq K}\langle \widehat{\mu_j}, \widehat{T_{\vec{\phi_r}}(\mu_k)}\eta_j\rangle\\
\nonumber &\lesssim& \sum_{|j-k| \leq K}\|\mu_{j}\|_2 \hskip.1cm \|\widehat{T_{\vec{\phi_r}}(\mu_k)}\eta_j\|_2\\
\end{eqnarray} where we use the Cauchy-Schwartz inequality.  Since $\mu$ is an Ahlfors-David regular measure on a set of Hausdorff dimension $\alpha$, that 
\begin{equation} \label{wolffnotes} \| \mu_j \|_2 \lesssim 2^{\frac{j(d-\alpha)}{2}}. \end{equation} 

Indeed, 
$$ {||\mu_j||}_2^2=\int {|\widehat{\mu}(\xi)|}^2 \eta(2^{-j} \xi) d\xi$$
$$=\int \int \int e^{2 \pi i (x-y) \cdot \xi} \eta(2^{-j} \xi) d\xi d\mu(x) d\mu(y)$$ 
$$=2^{dj} \int \int \widehat{\eta}(2^j(x-y)) d\mu(x) d\mu(y).$$ 

The absolute value of this quantity is bounded, for every $N>0$, by 
$$ C_N 2^{dj} \int \int {(1+2^j|x-y|)}^{-N} d\mu(x) d\mu(y)$$
$$=C_N 2^{dj} \int \int_{|x-y| \leq 2^{-j}} {(1+2^j|x-y|)}^{-N} d\mu(x) d\mu(y)$$
$$+C_N 2^{dj} \sum_{l=0}^{\infty} \int \int_{2^l \leq 2^j|x-y| \leq 2^{l+1}} {(1+2^j|x-y|)}^{-N} d\mu(x) d\mu(y)$$
$$=I+II.$$ 

By the Ahlfors-David property, 
$$ I \lesssim C_N 2^{dj} 2^{-j \alpha}.$$ 

Since $\mu$ is compactly supported, there exists $M>0$ such that 
$$ II=C_N 2^{dj} \sum_{l=0}^{j+M} \int \int_{2^l \leq 2^j|x-y| \leq 2^{l+1}} {(1+2^j|x-y|)}^{-N} d\mu(x) d\mu(y).$$

This expression is 
$$ \lesssim C_N 2^{dj} \sum_{l=0}^{j+M} 2^{-j \alpha} 2^{l \alpha} 2^{-lN} \lesssim C_N 2^{j(d-\alpha)}.$$ 

It follows that $I+II \lesssim 2^{j(d-\alpha)}$ and (\ref{wolffnotes}) is established. 

\vskip.125in 

We also have that
\begin{equation} \label{microlocal} \|\widehat{T_{\vec{\phi_r}}(\mu_k)}\eta_j\|_2 \lesssim 2^{-ks} 2^{\frac{k(d-\alpha)}{2}} \end{equation} by the mapping properties of the operator $T_{\vec{\phi_r}}$ in the regime of $|j-k|<K$.

The following lemma is a variant of a calculation in \cite{IJL10}. We will use it to get a bound on the second sum.

\begin{lemma} \label{fardiag}
For any $M>2d + m + 1$ there exists a constant $C_M>0$ such that for all indices $j,k$ with $|j-k|>K$ with $K$ large enough,

$$\langle T_{\vec{\phi_r}}\mu_j,\mu_k \rangle \leq C_M 2^{-M \max \{j,k\}}. $$ 
\end{lemma}

\noindent To prove the lemma, for simplicity, we replace $T_{\vec{\phi_r}}$ by $T$ and write 

$$T\mu_k(x)= \int_{\{y: \vec{\phi}(x,y)=r\}} \psi(x,y) \mu_k(y)d\sigma_{x,r}(y),$$
where $d\sigma_{x,r}$ is the Lebesgue measure on the set 
$\{y: \vec{\phi}(x,y)=r\}$.
\noindent It follows from our upcoming arguments that as long as $t_l \leq r_l \leq t_l + \varepsilon$ , the estimates hold uniformly in $r$.

As $\phi$ satisifies the property that $\{\nabla_y \phi_l(x,y)\}_{l=1}^m$ are linearly independent on a relatively open, bounded subset of $\{y: \phi(x,y)=t\}$ from (\ref{fibration}), we can assume that $|\sum_l \nabla_y \phi_l(x,y)| \approx 1$ on this set by making the support of $\psi$ small enough.  Next, we use an approximation argument on $T$ by letting 
\begin{equation}
T_n\mu_k(x)=n^m \int_{\mathbb{R}^d} \psi(x,y) \Pi_l \chi_l(n(\phi_l(x,y)-r_l)) \mu_k(y)dy
\end{equation} where $\{\chi_l\}_{l=1}^m$ is a family of smooth cutoffs supported near $0$ and equal to 1 near 0.  It is shown in \cite{GelShi66} that 
\begin{equation}
n^m\Pi_l \chi_l(n(\phi_l(x,y)-r_l)) dy
\end{equation} converges to the measure that appears in $T_{\vec{\phi}_r}$ as $n \rightarrow \infty$.  Therefore, proving the estimate in the case where $T_{\vec{\phi}_r}$ is replaced by $T_n$  is sufficient by convergence theorems found in \cite{Fol84} which in turn shows the uniformity in $r$.  We will drop the domains of integration in the upcoming calculations for brevity.

By Fourier inversion, we have 
$$T_n\mu(x)= \int e^{i y \cdot \xi} e^{is \cdot (\vec{\phi}(x,y)-r)}  \psi(x,y) \Pi_l\widehat{\chi_l}(n^{-1}s_l)  \widehat{\mu}(\xi) d\xi ds dy $$ 
\noindent and therefore
\begin{equation}\label{FourierTransT}\widehat{T_n \mu}(\eta)= \int  e^{-i x \cdot \eta} e^{i y \cdot \xi} e^{is \cdot (\vec{\phi}(x,y)-r)}  \psi(x,y) \Pi_l\widehat{\chi_l}(n^{-1}s_l)  \widehat{\mu}(\xi)  dx dy ds d\xi.\end{equation}

\noindent Invoking the properties of the Fourier transform on $L^2$, we see that 
\begin{eqnarray}
\langle T_n\mu_j,\mu_k \rangle &=& \langle \widehat{T_n \mu_j}, \widehat{\mu_k}\rangle \\
\nonumber &=& \int  e^{-i x \cdot \eta} e^{i y \cdot \xi} e^{is \cdot (\vec{\phi}(x,y)-r)}  \psi(x,y) \Pi_l\widehat{\chi_l}(n^{-1}s_l) \widehat{\mu_j}(\xi) \widehat{\mu_k}(\eta) dx dy ds d\xi d\eta \\
 \label{negativedecay} &=& \int \widehat{\mu_j}(\xi) \widehat{\mu_k}(\eta) \Pi_l\widehat{\chi_l}(n^{-1}s_l) I_{jk}(\xi, \eta, s) d \eta d \xi ds 
\end{eqnarray}

\noindent where 

\begin{equation} \label{oscillatory}
I_{jk}(\xi, \eta, s) = \psi_0(2^{-j}|\xi|)\psi_0(2^{-k}|\eta|) \int e^{is \cdot (\vec{\phi}(x,y)-r)}  e^{i y \cdot \xi} e^{-i x \cdot \eta} \psi(x,y) dx dy
\end{equation}

\noindent and $\psi_0$ is smooth cutoff equal to 1 on $\{1 \leq |z| \leq 10 \}$ and vanishing in the ball of radius 1/2.  The justification of such cutoffs comes from the support of $\widehat{\mu_j}(\xi)$ and $\widehat{\mu_k}(\eta)$ and again that $\eta_j \approx \eta_j^2$.  We will show that 
\begin{equation}\label{ibpdecay}|I_{jk}(\xi, \eta, s)| \leq C_M 2^{-M max(j,k)}\end{equation} when $|j-k|>K$ for a large enough $K$.

Computing the critical points of the phase function in (\ref{oscillatory}), we see that 
\begin{equation}
 \sum_l |s|\tilde{s}_l \nabla_x\phi_l(x,y)=\eta  \text{ and } \sum_l |s|\tilde{s}_l \nabla_y \phi_l(x,y)=-\xi,\end{equation} where $s=|s|(\tilde{s}_1,...,\tilde{s}_m)$ and $(\tilde{s}_1,...,\tilde{s}_m) \in S^{m-1}$, the unit sphere. The compactness of the support of $\psi$ and the domain of the variable $(\tilde{s}_1,...,\tilde{s}_m)$ along with the linear independence condition from (\ref{fibration}) implies that 
\begin{equation}
 \left|\sum_l \tilde{s}_l \nabla_x\phi_l(x,y)\right| \approx \left|\sum_l \tilde{s}_l \nabla_y\phi_l(x,y)\right| \approx 1.
 \end{equation}  More precisely, the upper bound follows from smoothness and compact support. The lower bound follows from the fact that a continuous non-negative function achieves its minimum on a compact set. This minimum is not zero because of the linear independence condition (\ref{fibration}). 

It follows that 
\begin{equation} \label{climax} |\xi| \approx |\eta| \end{equation} when we are near the critical points in $(x,y)$. The support of the cutoffs $\psi_0$, when $|j-k|>K$, tell us that we are supported away from critical points in $(x,y)$  since (\ref{climax}) no longer holds.  This condition implies that for some $h$ or $h'$ in $\{1,2,...,d\}$, 
\begin{equation} \left(\sum_l s_l \frac{\partial \phi_l}{\partial x_h} - \eta_h \right) \neq 0 \ \text{or} \ \left(\sum_l s_l \frac{\partial \phi_l}{\partial y_{h'}} + \xi_{h'}\right) \neq 0.
\end{equation}

Without loss of generality, assume the former holds and that $k> j$.  It is immediate that $e^{-i x \cdot \eta} e^{is \cdot(\vec{\phi}(x,y)-r)}$ is an eigenfunction of the differential operator 
\begin{equation}
L = \frac{1}{i(\sum_l s_l \frac{\partial \phi_l}{\partial x_h} - \eta_h)} \frac{\partial}{\partial x_h};
\end{equation}  We integrate by parts in (\ref{oscillatory}) using this operator.  The expression that we get after performing this procedure $M > 2d + m + 1$ times is 
\begin{equation}
I(\xi, \eta, s) \lesssim \sup_{x,y} \left| \sum_l s_l \frac{\partial \phi_l}{\partial x_h} - \eta_h \right|^{-M}.
\end{equation}

\noindent Now, suppose that we are in the region $\{|s| << |\eta|\}$ (i.e $|s| \leq c|\eta|$ with a sufficiently large constant $c>0$).  Since $| \sum_l s_l \nabla_x \phi_l| \approx |s|$ it follows, after possibly changing our initial choice of $h$, that 
\begin{equation}
\left|\sum_l s_l \frac{\partial \phi_l}{\partial x_h} - \eta_h \right| \gtrsim \left| \hskip.1cm \left|\sum_l s_l \frac{\partial \phi_l}{\partial x_h}\right| - |\eta| \hskip.1cm \right| \approx |\eta|.
\end{equation} Similarly, if $\{|s| >> |\eta|\}$ then, again after possibly changing our initial choice of $h$,
\begin{equation}
\left|\sum_l s_l \frac{\partial \phi_l}{\partial x_h} - \eta_h \right| \gtrsim \left| \hskip.1cm \left|\sum_l s_l \frac{\partial \phi_l}{\partial x_h}\right| - |\eta| \hskip.1cm \right| \approx |s|.
\end{equation} 

\noindent In either region, 
\begin{equation}
|I_{jk}(\xi, \eta, s)| \lesssim \sup(|s|, |\eta|)^{-M} \lesssim 2^{-Mk} .
\end{equation}

Considering (\ref{negativedecay}), the integrand $(\Pi_l\widehat{\chi_l}(n^{-1}s_l)) I_{jk}(\xi, \eta, s)$ is integrable in $s$ as the first term is at most $1$ and $I_{jk}$ is bounded about by $|s|^{-M}$.  Performing the remaining integrations and keeping in mind the support properties of $\widehat{\mu_j}$ and $\widehat{\mu_k}$, it follows that 
\begin{equation} \label{2.9}
\sum_{|j-k| > K} \langle \mu_j, T_{\vec{\phi}_r}(\mu_k)\rangle \lesssim \sum_{|j-k| > K} C_M 2^{-(M-2d) max(j,k)} \lesssim 1.
\end{equation}
This completes the proof of Lemma \ref{fardiag}.

We are now ready to give the final step of the proof to Proposition \ref{mamaestprop}.  Since both sums in (\ref{2.6}) are bounded by 1, this implies that the left hand side of (\ref{2.3}) is bounded above by $\varepsilon^m$ after completing the integrations in (\ref{expanded}).

\vskip.25in 

\section{Proof of the Corollary (\ref{phongstein})} 
\label{corollarysection} 

In order to establish Corollary \ref{phongstein}, we need to prove that the estimate (\ref{sobolev}) holds with $s=\frac{d-1}{2}$ under the Phong-Stein condition (\ref{curvecond}). This follows, for example, from the main result in \cite{PhSt86}. See also \cite{Ho71} and \cite{So93} for a thorough description of related estimates. In this section, we place this estimate into the context of general Fourier integral operator theory for the sake of clarity. 

The Radon transform, which sends a function to its averages on a given family of submanifolds, has appeared frequently in many areas of analysis and geometry.  Its appearance, for example, in the study of the $\overline{\partial}$-Neumann problem and integral geometry \cite{PhSt86} brought microlocal analysis past its initial uses in the analysis of parametrices and the propogation of singularities.  The condition (\ref{curvecond}) can be viewed as a nondegeneracy assumption when taking  H\"ormander's viewpoint of Fourier integral operators \cite{Ho71}. 

Let us consider the integral operator 

\begin{equation} 
\mathfrak{T}f(x)=\int K(x,y)f(y) dy
\end{equation}
 
\noindent where $f \in L^2(\mathbb{R}^d)$ and 

\begin{equation} \label{kernel}
K(x,y) = \int e^{i \Psi(x,y,\theta)} a(x,y,\theta) d \theta
\end{equation} 

\noindent for $\Psi \in C^{\infty}(\mathbb{R}^d \times \mathbb{R}^d \times \mathbb{R}^N)$ and $a \in C_0^{\infty}(\mathbb{R}^d \times \mathbb{R}^d \times \mathbb{R}^N)$.  

The $L^2$ mapping properites of $\mathfrak{T}$ are determined by the geometric properties of the canonical relation 

\begin{equation} 
C = \{(x, \nabla_x \Psi, y, \nabla_y \Psi) \hskip.1cm : \hskip.1cm \nabla_{\theta} \Psi = 0 \}
\end{equation}

\noindent which is a subset of $T^*(\mathbb{R}^d \times \mathbb{R}^d)$.  If the set $\{\nabla_{x,y,\theta} (\partial \Psi / \partial \theta_i) \}_{i=1}^d$ is linearly independent on $\{\nabla_{\theta} \Psi = 0\}$, then $C$ is an immersed submanifold.  Moreover, these conditions put $K$ into the general framework of Fourier integral distributions \cite{Ho71}.  We now call the operator $\mathfrak{T}$ a Fourier integral operator associated to $C$.  When $T_{\phi_t}$ is viewed as a Fourier integral operator, 
(\ref{kernel}) becomes 

\begin{equation} \label{kernelT}
\int e^{i (\phi(x,y)-t)\tau} a(x,y,\tau) d \tau 
\end{equation}

\noindent for $\tau \in \mathbb{R}$. 

The best possible situation for $L^2$ estimates for $\mathfrak{T}$ comes when $C$ is locally the graph of a canonical transformation \cite{Ho71}.  This is equivalent to 

\begin{equation}
\begin{pmatrix}
 Hess_{xy}\Psi & Hess_{x \theta}\Psi \\
 Hess_{\theta y}\Psi & Hess_{\theta \theta}\Psi
\end{pmatrix}
\neq 0
\end{equation}

\noindent where $Hess_{z',z''}\Psi$ is the mixed Hessian of $\Psi$ in the variables $z'$ and $z''$.  The resulting $L^2$ estimate for $\mathfrak{T}$ is $L^2(\mathbb{R}^d) \rightarrow L^2_{\frac{d-N}{2}}(\mathbb{R}^d)$.  Computing determinant in the case of (\ref{kernelT}), we get the Monge-Ampere determinant that appears in the definition of the Phong-Stein rotational curvature condition.  Hence, (\ref{curvecond}) guarantees that $T_{\phi_t}$ has as its canonical relation a local canonical graph and is smoothing of order $\frac{d-1}{2}$ on $L^2({\Bbb R}^d)$.  

By (\ref{dimensionrestriction}) and the fact that $s = \frac{d-1}{2}$, it follows that $\alpha > \frac{d+1}{2}$ and the proof of Corollary \ref{phongstein} is complete. 

\vskip.25in 

\section{Sharpness of results} 
\label{finalsection} 

\subsection{Proof of Theorem \ref{sharp}} 

We use a construction very closely related to the one in \cite{IS10}. However, we will start with a variant of the incidence example due to Pavel Valtr \cite{BMP00, V05} as an exercise for our intuition.  Let us note that a similar object can be found in \cite{K77} in a slightly different context. Let 
\begin{equation} \label{valt1}
 P_n=\left\{\left(\frac{i_1}{n}, \frac{i_2}{n}, \dots, \frac{i_{d-1}}{n}, \frac{i_d}{n^2}\right):0\leq i_j\leq n-1,\text{ for } 1\leq j\leq d-1, \text{ and } 1\leq i_d \leq n^2\right\}.
 \end{equation}
Notice that in each of the first $d-1$ coordinates, there are $n$ evenly distributed points, but in the last dimension, there are $n^2$ evenly distributed points. Now, let 
\begin{equation} \label{valt2}
 H=\{(t_1,t_2, \dots, t_{d-1}, t_1^2+\dots+t_{d-1}^2) \in {\Bbb R}^d: t_l \in {\Bbb R}\}
 \end{equation} and define 
\begin{equation} \label{valt3}
 L_H=\{H+p, p \in P_n\}.
 \end{equation}  Note that $L_H$ is a collection of shifted paraboloids. 

\begin{figure}
\centering
\includegraphics[scale=1]{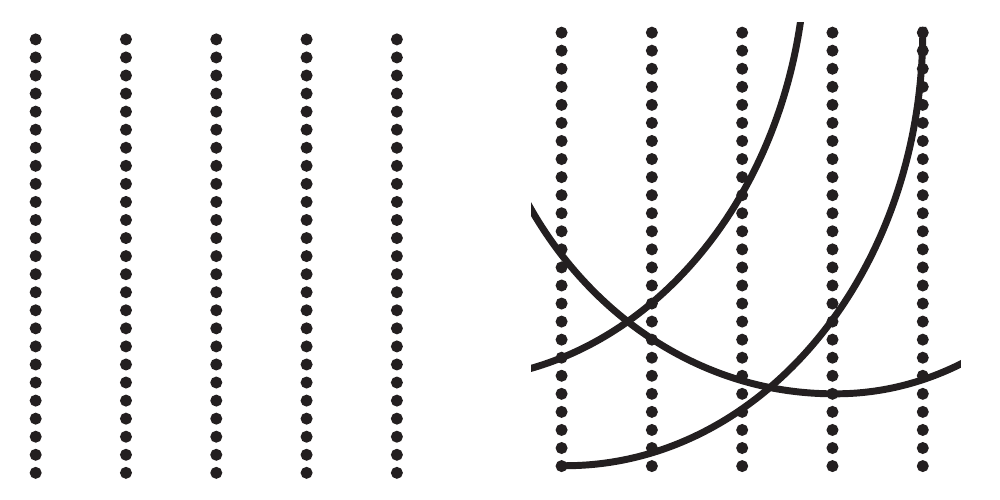}
\caption{On the left, we see a picture of the set $P_5$, on the right, we see it again with a few parabolic arcs, which intersect a point in each column.}
\label{ValtrFig1}
\end{figure}

Let $N=n^{d+1}$. By construction, $\# P_n= \# L_H = N$. Also by construction, each element of $L_H$ is incident to about $n^{d-1} \approx N^{\frac{d-1}{d+1}}$ elements of $P_n$. Thus the total number of incidences between $P_n$ and $L_H$ is 
$$\approx N^{1+\frac{d-1}{d+1}}=N^{\frac{2d}{d+1}}=N^{2-\frac{2}{d+1}}.$$

We are now ready to define the convex body $B$. With the combinatorial construction in (\ref{valt1}), (\ref{valt3}), and (\ref{valt3}) in hand, we flip the paraboloid upside down and glue it to another copy. Explicitly, let
$$B_U=\left\lbrace(x_1, x_2, \dots, x_d) \in \mathbb{R}^d : x_i \in [-1,1], \text{ for } 1 \leq i \leq d-1, \text{ and } x_d = 1-\left( x_1^2+x_2^2+ \dots + x_{d-1}^2 \right) \right\rbrace,$$
and
$$B_L=\left\lbrace(x_1, x_2, \dots, x_d) \in \mathbb{R}^d : x_i \in [-1,1], \text{ for } 1 \leq i \leq d-1, \text{ and } x_d = -1+ x_1^2+x_2^2+ \dots + x_{d-1}^2 \right\rbrace.$$
Now, let
$$B'=\left(B_U \cap \left\lbrace (x_1, x_2, \dots, x_d) \in \mathbb{R}^d :x_d \geq 0\right\rbrace \right) \cup \left( B_L \cap \left\lbrace (x_1, x_2, \dots, x_d) \in \mathbb{R}^d :x_d \leq 0 \right\rbrace \right).$$
Finally, define $B$ to be the convex body $B'$, with the ridge at the transition between $B_U$ and $B_L$ smoothed.

We now transition to show how the upper Minkowski dimension of a particular $(S_t^{\phi}(E))$ for an $E$ given below gives the desired sharpness result.

Let $\{q_i\}_{i\in \mathbb{N}}$ be a sequence of positive integers such that $q_{i+1}=q_i^i$ and $q_1=2$. Let $E_i$ be the $q_i^{-\frac{d}{s}}$ neighborhood of the set 
$$ q_i^{-1} \left\{x \in {\Bbb Z}^d: 0 \leq x_j \leq q_i^{\frac{d}{d+1}}, 1 \leq j \leq d-1, 0 \leq x_d \leq q_i^{\frac{2d}{d+1}} \right\},$$ where ${\Bbb Z}^d$ denotes the standard integer lattice.  We set $E=\cap_i E_i$.  See \cite{Falc86}, Chapter 8 for a proof that $E=\cap_i E_i$ is an Ahlfors-David regular set of Hausdorff dimension $s$.  

Define $\phi(x,y):={||x-y||}_B$, where ${|| \cdot ||}_B$ is the norm induced by a symmetric convex body $B$ with a smooth boundary and non-vanishing Gaussian curvature that we now construct. A fairly straightforward calculation shows that as long as $\partial B$ is smooth and has everywhere non-vanishing Gaussian curvature, then $\phi(x,y)={||x-y||}_B$ satisfies the Phong-Stein rotational curvature hypothesis away from $t=0$. 

By the definition of Minkowski dimension \cite{W04}, our problem of calculating $\overline{dim}_{{\mathcal M}}(S_t^{\phi}(E))$ reduces to counting the number of $\delta$-balls (or boxes) needed to cover $S_t^{\phi}(E)$. We will see that the number of $\delta_i$-boxes needed to cover $E_i$ is the same as the number needed to cover $E$. 
 
Set $\delta_i = q_i^{-\frac{d}{s}}$ and partition the set $[0,1]^d$ into boxes of sidelength $\sim \delta_i$.
The set $E$ being Ahlfors-David regular implies that $\mu(B(x, \delta_i))\sim \delta_i^s$ whenever $x\in E$.  Hence there are $\sim \delta_i^{-s}$ boxes of the $\delta_i^{-d}$ boxes from the partioning of $[0,1]^d$, each of which contains a non-isolated point of $E$.  The definition of $E = \cap_i E_i$ shows that the number of $\delta_i$-boxes needed to cover $E_i$ is the same as the number of $\delta_i$ boxes needed to cover $E$. Our problem of calculating $\overline{dim}_{{\mathcal M}}(S_t^{\phi}(E))$ now reduces to the problem of calculating the number of balls of radius $q_i^{-\frac{d}{s}}$ needed to cover $S_t^{\phi}(E_i)$. 

Given $x \in E_i$ and keeping mind the definition of $E$, the number of balls of radius $q_i^{-\frac{d}{s}}$ needed to cover $\{y: {||x-y||}_B=1 \}$ is $\gtrsim q_i^{\frac{d(d-1)}{d+1}}$. It follows that the number of balls of radius $q_i^{-\frac{d}{s}}$ needed to cover $S_t^{\phi}(E_i)$ and, consequently, $S_t^{\phi}(E)$, is 
$$ \gtrsim q_i^d \cdot q_i^{\frac{d(d-1)}{d+1}}=q_i^{\frac{2d^2}{d+1}}={(q_i^{-\frac{d}{s}})}^{-\frac{2ds}{d+1}}.$$

Therefore the upper Minkowski dimension of $S_t^{\phi}(E)$ is at least $\frac{2ds}{d+1}$. This number is greater than $2s-1$ if $s<\frac{d+1}{2}$. 

\subsection{Sharpness of the method}

It is known that the estimate (\ref{mama}), which is at the core of our method, is essentially sharp in the case $m=1$. See \cite{Mat85, Mat87} for the proof of this fact in the $d=2$ case and \cite{IS10} for the construction in dimensions $d \geq 3$.  Mattila's example is for the function 
$\phi(x,y)=|x-y|$, where $| \cdot |$ is the Euclidean distance. The construction in \cite{IS10} is for the convex body $B$ used in the proof of Theorem \ref{sharp} above. We now give a proof of sharpness of the exponent 
$\frac{3}{2}$ in (\ref{mama}) in the case $d=2$ and $\phi(x,y)=x \cdot y$ as this result appears to be new. We believe that when $m=1$, the exponent $\frac{d+1}{2}$ is sharp in all dimensions for any function $\phi$ satisfying the Phong-Stein rotational curvature condition, but are unable to prove this at the moment. 

Consider the estimate 
\begin{equation} \label{dot} \mu \times \mu \{(x,y): 1 \leq x\cdot y \leq 1+\epsilon \}\lesssim \epsilon. \end{equation} We will modify the referenced examples to show that for no $s<\frac{3}{2}$ does 
\begin{equation}
I_s(\mu)=\int \int {|x-y|}^{-s} d\mu(x)d\mu(y)<\infty
\end{equation} imply (\ref{dot}).

We denote by $\mathcal{C}_\alpha$, a Cantor set of Hausdorff dimension $0<\alpha<1$.  Let $F=\left(\mathcal{C}_\alpha\cap [\frac{1}{2},1]\right)\cup \left(\mathcal{C}_\alpha^{-1}\cap[\frac{1}{2},1]\right)$.  Now define $\mathcal{M}_2(\alpha) = \left\{r\omega : r\in F, \omega \in S^1 \right\}  \subset \mathbb{R}^2.$  Set $\mu = \mathcal{H}^\alpha_{|F}  \times \mathcal{L}^1_{|[0,1]}$ where $\mathcal{H}^\alpha$ is the $\alpha$-dimensional Hausdorff measure.  This measure is meant to be a fractal analog of that from polar coordinates in the plane.

Pick a point $x=r\omega \in \mathcal{M}_2(\alpha)$. Notice that if $r\in F$, then so is $\frac{1}{r}$. Notice that $\left\{y:1\leq x\cdot y\leq 1+\epsilon\right\}$ is contained in a strip formed by the two lines which are both perpendicular to the vector $x$ and pass through the points $\frac{1}{r}\omega$ and $\frac{}{}\omega$ respectively. 
  
We argue that within an $\epsilon$-annulus we can fit a rectangle of width $\sim \epsilon$ and length $\sim\sqrt{\epsilon}$. Similarly, an $\epsilon$-strip fits an annulus of width $\sim \epsilon$ and inner arc-length $\sim\sqrt{\epsilon}$.  This rectangle intersects $\mathcal{M}_2(\alpha)$ for $x$ in a set of positive $\mu$ measure. The measure of this intersection is $\sim \epsilon^{1/2+\alpha}$. It follows that
\begin{equation}
 \mu\left\{y:1\leq x\cdot y\leq1+\epsilon\right\} \gtrsim \epsilon^{\alpha+1/2} 
\end{equation} for $x$ in a set of positive $\mu$ measure and consequently
\begin{eqnarray}
 \mu \times \mu \{(x,y): 1 \leq x\cdot y \leq 1+\epsilon \} &=& \int  \mu\left\{y:1\leq x\cdot y \leq1+\epsilon\right\} d\mu(x)\\
\nonumber && \gtrsim \epsilon^{\alpha+1/2}.
\end{eqnarray}  Hence,  
\begin{equation}
 \mu \times \mu \{(x,y): 1 \leq x\cdot y \leq 1+\epsilon \} \lesssim \epsilon
 \end{equation} only for 
\begin{equation}
\epsilon^{\alpha+\frac{1}{2}} \lesssim \epsilon,
\end{equation} which can only hold when $\alpha \ge \frac{1}{2}.$  Thus the estimate (\ref{dot}) does not in general hold for sets with Hausdorff dimension less than $\frac{3}{2}$.

\end{document}